\documentclass[11pt]{article}
\usepackage{amsthm,amsmath,latexsym,amssymb}

\newcommand{\bP}{{\rm |\kern-.15em P}}
\newcommand{\Q}{\kern.3em\rule{.07em}{.65em}\kern-.3em{\rm Q}}
\newcommand{\R}{{\rm I\kern-.15em R}}
\newcommand{\D}{{\rm |\kern-.15em D}}
\newcommand{\h}{{\rm |\kern-.15em H}}
\newcommand{\C}{\kern.3em\rule{.07em}{.65em}\kern-.3em{\rm C}}
\newcommand{\T}{{\rm T\kern-.35em T}}

\theoremstyle{plain}
\newtheorem{theorem}{Theorem}[section]

\newtheorem{corollary}[theorem]{Corollary}

\theoremstyle{definition}
\newtheorem{definition}[theorem]{Definition}

\theoremstyle{remark}
\newtheorem{remark}[theorem]{Remark}

\begin{document}
\title{Fibers of Polynomial Mappings Over $\mathbb{R}^n$}
\author{Ronen Peretz}
 
\maketitle

\begin{abstract}
We prove results on fibers of polynomial mappings $\mathbb{R}^n\rightarrow\mathbb{R}^n$ and
deduce when such mappings are surjective under certain conditions.
\end{abstract}

\section{The results}

\begin{definition}
Let $p(X_1,\ldots,X_n)\in\mathbb{R}[X_1,\ldots,X_n]$. We denote by $\overline{p}(X_1,\ldots,X_n)$ the leading homogeneous 
component of $p(X_1,\ldots,X_n)$ with respect to the standard grading, $\deg X_j=1$ for $1\le j\le n$.
\end{definition}

\begin{theorem}
Let $f:\,\mathbb{R}^n\rightarrow\mathbb{R}^n$, $f(X_1,\ldots,X_n)=(p_1(X_1,\ldots,X_n),\ldots,p_n(X_1,\ldots,X_n))$ be a polynomial mapping
(i.e. $(p_1,\ldots,p_n)\in\mathbb{R}[X_1,\ldots,X_n]^n$). Let $g_{ij}(X_1,\ldots,X_n)\in\mathbb{R}[X_1,\ldots,X_n]$, $i,j=1,\ldots,n$.
Let $(\alpha_1,\ldots,\alpha_n)\in (\mathbb{Z}^+)^n$. We assume that the following $2$ conditions hold true: \\
(i) The determinant $\det\left(g_{ij}(X_1,\ldots,X_n)\right)_{i,j=1,\ldots,n}$ never vanishes in $\mathbb{R}^n$. \\
(ii) The following system of $n$ equations in $n$ unknowns is such that the degree of each of the equations is an odd number:
\begin{equation}
\label{l1}
\sum_{j=1}^n\left(p_j(X_1,\ldots,X_n)\right)^{\alpha_j}g_{ij}(X_1,\ldots,X_n)=0,\,\,\,i=1,\ldots,n.
\end{equation}
Then the following $2$  assertions are true: \\
(a) If the induced homogeneous system of the system (\ref{l1}):
\begin{equation}
\label{l2}
\overline{\sum_{j=1}^n\left(p_j(X_1,\ldots,X_n)\right)^{\alpha_j}g_{ij}(X_1,\ldots,X_n)}=0,\,\,\,i=1,\ldots,n,
\end{equation}
has only the zero solution $(X_1,\ldots,X_n)=(0,\ldots,0)$ over $\mathbb{R}$, then $f(\mathbb{R}^n)=\mathbb{R}^n$. \\
(b)If the induced homogeneous system of the system in equation (\ref{l1}), i.e. the system (\ref{l2}) has only the zero
solution over $\mathbb{C}$, then $\forall\,(a_1,\ldots,a_2)\in\mathbb{R}^n$, either $|f^{-1}(a_1,\ldots,a_n)|=\infty$
over $\mathbb{C}$ under the extra assumption that $\det\left(g_{ij}(Z_1,\ldots,Z_n)\right)_{i,j=1,\ldots,n}\in\mathbb{R}^{\times}$,
or there exists an integer $k=k(a_1,\ldots,a_n)\ge 0$ such that $|f^{-1}(a_1,\ldots,a_n)|=2k+1$ over $\mathbb{R}$.
\end{theorem}
\noindent
{\bf Proof.} \\
(a) Let $(a_1,\ldots,a_n)\in\mathbb{R}^n$. We will prove that $(a_1,\ldots,a_n)\in f(\mathbb{R}^n)$. we consider the following system
of equations:
$$
X_{n+1}^{d_i}\sum_{j=1}^n\left(p_j\left(\frac{X_1}{X_{n+1}},\ldots,\frac{X_n}{X_{n+1}}\right)-a_j\right)^{\alpha_j}
g_{ij}\left(\frac{X_1}{X_{n+1}},\ldots,\frac{X_n}{X_{n+1}}\right)=0,
$$
\begin{equation}
\label{l3}
\end{equation}
$$
{\rm where}\,\,\,\,\,d_i=\deg\left(\sum_{j=1}^n \left(p_j\right)^{\alpha_j}g_{ij}\right),\,\,\,\,\,i=1,\ldots,n.
$$
This is a system of $n$homogeneous real polynomial equations in the $n+1$ unknowns $X_1,\ldots,X_n,X_{n+1}$, and by condition (ii)
the degrees $d_i$, $i=1,\ldots,n$ of all of these equations are odd integers. By well known facts on varieties over $\mathbb{R}$
(see \cite{shaf}), it follows that the system (\ref{l3}) has a non-zero real solution $(X_1,\ldots,X_n,X_{n+1})=
(X_1^0,\ldots,X_n^0,X_{n+1}^0)$. We must have $X_{n+1}^0\ne 0$, for otherwise $(X_1^0,\ldots,X_n^0)\ne (0,\ldots,0)$ and 
$(X_1^0,\ldots,X_n^0)$ is a solution of (\ref{l1}),i.e. (\ref{l2}). This contradicts the assumption of the theorem in part (a).
Thus we get from equation (\ref{l3}):
$$
\sum_{j=1}^n\left(p_j\left(\frac{X_1^0}{X_{n+1}^0},\ldots,\frac{X_n^0}{X_{n+1}^0}\right)-a_j\right)^{\alpha_j}
g_{ij}\left(\frac{X_1^0}{X_{n+1}^0},\ldots,\frac{X_n^0}{X_{n+1}^0}\right)=0,
$$
\begin{equation}
\label{l4}
{\rm for}\,\,\,\,\,i=1,\ldots,n.
\end{equation}
By condition (i) of our theorem, this implies that
$$
f\left(\frac{X_1^0}{X_{n+1}^0},\ldots,\frac{X_n^0}{X_{n+1}^0}\right)=(a_1,\ldots,a_n).
$$
(b) Let us consider the system (\ref{l3}) over $\mathbb{C}$. By the Bezout Theorem (see \cite{shaf}), either the system (\ref{l3})
has infinitely many solutions over $\mathbb{C}$, or it has exactly
$$
\prod_{i=1}^n\deg\left(\sum_{j=1}^n(p_j)^{\alpha_j}g_{ij}\right)
$$
solutions over $\mathbb{C}$, counting multiplicities and not counting the zero solution. In the case we have infinitely many solutions
over $\mathbb{C}$, we must have for each such a solution $(Z_1^0,\ldots,Z_n^0,Z_{n+1}^0)$ that $Z_{n+1}^0\ne 0$, for by the 
assumption in part (b) of our theorem, the induced homogeneous system (\ref{l2}), of the system (\ref{l1}) has only the zero
solution over $\mathbb{C}$. Since we also assume in this case that $\det(g_{ij}(Z_1,\ldots,Z_n))_{i,j=1,\ldots,n}\in\mathbb{C}^{\times}$
it follows as before by equation (\ref{l4}) that the fiber over $\mathbb{C}$, $f^{-1}(a_1,\ldots,a_n)$ contains infinitely many points:
$$
\left(\frac{Z_1^0}{Z_{n+1}^0},\ldots,\frac{Z_n^0}{Z_{n+1}^0}\right).
$$
In the second case, in which we have exactly 
$$
\prod_{i=1}^n\deg\left(\sum_{j=1}^n(p_j)^{\alpha_j}g_{ij}\right)
$$
solutions over $\mathbb{C}$, noting that by condition (ii) this number is an odd integer and that non-real solutions
$(Z_1^0,\ldots,Z_n^0,Z_{n+1}^0)$ come in conjugate pairs, we deduce that the fiber over $\mathbb{R}$, $f^{-1}(a_1,\ldots,a_n)$,
contains an odd number of points. $\qed $ \\
\\
\begin{corollary}
Let the polynomial mapping $f:\,\mathbb{R}^n\rightarrow\mathbb{R}^n$, be given by $f(X_1,\ldots,X_n)=(p_1(X_1,\ldots,X_n),\ldots,
p_n(X_1,\ldots,X_n))$. Let $g_{ij}(X_1,\ldots,X_n)\in\mathbb{R}[X_1,\ldots,X_n]$ for $i,j=1,\ldots,n$. 
Let $(\alpha_1,\ldots,\alpha_n)\in \left(\mathbb{Z}^+\right)^n$.
We assume that the following two conditions hold true: \\
(i) The determinant $\det\left(g_{ij}(X_1,\ldots,X_n)\right)_{i,j=1,\ldots,n}$ never vanishes in $\mathbb{R}^n$. \\
(ii) For each $i=1,\ldots,n$ the set $\{\alpha_j\deg p_j+\deg g_{ij}\,|\,j=1,\ldots,n\}$ contains a unique maximal element
$\alpha_{j(i)}\deg p_{j(i)}+\deg g_{ij(i)}$, which is an odd integer. We agree that $\deg 0=-\infty$. \\
Let us consider the following homogeneous system:
\begin{equation}
\label{l5}
\overline{p}_{j(i)}\overline{g}_{ij(i)}=0,\,\,\,\,\,i=1,\ldots,n.
\end{equation}
Then the following two assertions are true: \\
(a) If the system (\ref{l5}) has only the zero solution over $\mathbb{R}$, then $f(\mathbb{R}^n)=\mathbb{R}^n$. \\
(b) If the system (\ref{l5}) has only the zero solution over $\mathbb{C}$, then for any $(a_1,\ldots,a_n)\in\mathbb{R}^n$
either $|f^{-1}(a_1,\ldots,a_n)|=\infty$ over $\mathbb{C}$, provided that also the following assumption holds true,
$\det g_{ij}(Z_1,\ldots,Z_n))_{i,j=1,\ldots,n}\in\mathbb{R}^{\times}$, or that there exists an integer $k=k(a_1,\ldots,a_n)\ge 0$
such that $|f^{-1}(a_1,\ldots,a_n)|=2k+1$ over $\mathbb{R}$.
\end{corollary}
\noindent
{\bf Proof.} \\
This is a special case of Theorem 1.2, where the system (\ref{l5}) is precisely the system (\ref{l2}) because of the maximality and 
the uniqueness of $\alpha_{j(i)}\deg p_{j(i)}+\deg g_{ij(i)}$ among the elements of the set 
$\{\alpha_j\deg p_j+\deg g_{ij}\,|\,j=1,\ldots,n\}$. $\qed $ \\
\\
\begin{corollary}
Let the polynomial mapping $f:\,\mathbb{R}^n\rightarrow\mathbb{R}^n$, be given by $f(X_1,\ldots,X_n)=(p_1(X_1,\ldots,X_n),\ldots,
p_n(X_1,\ldots,X_n))$. Suppose that the product $(\deg p_1)\cdot\ldots\cdot (\deg p_n)$ is an odd integer. Then the following two assertions
are true: \\
(a) If $|\overline{f}^{-1}(0,\ldots,0)|=1$ over $\mathbb{R}$, then $f(\mathbb{R}^n)=\mathbb{R}^n$. \\
(b) If $|\overline{f}^{-1}(0,\ldots,0)|=1$ over $\mathbb{C}$, then $\forall\,(a_1,\ldots,a_n)\in\mathbb{R}^n$ either the fiber size 
$|f^{-1}(a_1,\ldots,a_n)|=\infty$ over $\mathbb{C}$, or there exists an integer $k=k(a_1,\ldots,a_n)\ge 0$ such that
$|f^{-1}(a_1,\ldots,a_n)|=2k+1$ over $\mathbb{R}$.
\end{corollary}
\noindent
{\bf Proof.} \\
This follows by Corollary 1.3, where $(\alpha_1,\ldots,\alpha_n)=(1,\ldots,1)$ and where $g_{ij}=\delta_{ij}$, $i,j=1,\ldots,n$ because 
the system (\ref{l5}) becomes $\overline{p}_j=0$, $j=1,\ldots,n$ which has the solution set $\overline{f}^{-1}(0,\ldots,0)$. $\qed $ \\
\\
\begin{remark}
We note that if in Corollary 1.4 we have $\deg p_j=1$, $j=1,\ldots,n$, i.e. if all the $p_j=\overline{p}_j$ are linear forms
then we get the well known fact from linear algebra. Namely, if $A\overline{X}=\overline{0}$ is an $n\times n$ linear homogeneous
system that has only the trivial solution, then $A\overline{X}=\overline{b}$ is consistent $\forall\,\overline{b}\in\mathbb{R}^n$.
\end{remark}

\begin{remark}
If for $j=1,\ldots,n$, $bj\ge 0$ is an integer and if we have 
$$
p_j(X_1,\ldots,X_n)=\sum_{i=1}^n a_{ij}X_i^{2b_j+1}+{\rm elements}\,{\rm of}\,{\rm degrees}<2b_j+1.
$$
Then the polynomial mapping $f(X_1,\ldots,X_n)=(p_1(X_1,\ldots,X_n),\ldots,p_n(X_1,\ldots,X_n))$ is a surjective mapping, i.e.
$f(\mathbb{R}^n)=\mathbb{R}^n$, provided that the only solution of the following system:
$$
\sum_{i=1}^n a_{ij}X_i^{2b_j+1}=0,\,\,\,\,\,j=1,\ldots,n,
$$
is the trivial solution: $X_1=\ldots=X_n=0$. \\
For in this case the above system is the system (\ref{l5}) of Corollary 1.4 ($g_{ij}=\delta_{ij}$). For example, this is the case
for the equal-degree case $b_1=\ldots=b_n=b$ provided that $\det(a_{ij})_{i,j=1,\ldots,n}\ne 0$. Another example is the following: \\
we pick $4$ non-zero real numbers, $a,\,b,\,c$ and $d$ such that ${\rm sgn}(ad)=-{\rm sgn}(bc)$. Then any mapping of the form:
$$
f:\,\mathbb{R}^2\rightarrow\mathbb{R}^2,\,\,\,\,\,f(X,Y)=(aX^{2k+1}+bY^{2k+1}+\ldots,cX^{2j+1}+dY^{2j+1}+\ldots),
$$
is a surjective mapping. For the system (\ref{l5}) is:
$$
\left\{\begin{array}{ccc} aX^{2k+1}+bY^{2k+1} & = & 0 \\ cX^{2j+1}+dY^{2j+1} & = & 0 \end{array}\right..
$$
If $k\le j$ then the system can be written as follows:
$$
\left\{\begin{array}{lll} aX^{2k+1}+bY^{2k+1} & = & 0 \\ (cX^{2(j-k)})X^{2k+1}+(dY^{2(j-k)})Y^{2k+1} & = & 0 \end{array}\right..
$$
We view this as a linear homogeneous system in the unknowns $X^{2k+1}$ and $Y^{2k+1}$. Then the coefficients matrix is:
$$
\left(\begin{array}{cc}a & b \\ cX^{2(j-k)} & dY^{2(j-k)}\end{array}\right).
$$
The determinant of this matrix  is $(ad)Y^{(2(j-k)}-(bc)X^{2(j-k)}$ and this can not be $0$ because
of the assumption ${\rm sgn}(ad)=-{\rm sgn}(bc)$, unless $j>k$ and $X=Y=0$. In the other cases the only solution is, again,
$X=Y=0$.
\end{remark}

\begin{theorem}
Let $g_{ij}(X_1,\ldots,X_n)\in\mathbb{R}[X_1,\ldots,X_n]$ for $i,j=1,\ldots,n$ satisfy the condition that $\det(g_{ij}(X_1,\ldots,X_n))_{i,j=1,\ldots,n}$
never vanishes in $\mathbb{R}^n$. Then for any $j_0$, $1\le j_0\le n$, such that the degrees $\deg g_{ij_0}$, $i=1,\ldots,n$ are all odd integers 
the system:
\begin{equation}
\label{l6}
\overline{g}_{ij_0}(X_1,\ldots,X_n)=0,\,\,\,\,\,i=1,\ldots,n,
\end{equation}
has non-zero real solutions.
\end{theorem}
\noindent
{\bf Proof.} \\
Let $j_0$ be such that the degrees $\deg g_{ij_0}$, $i=1,\ldots,n$, are all odd integers. In Corollary 1.3 we take the following: 
$$
f:\,\mathbb{R}^n\rightarrow\mathbb{R}^n,\,\,\,\,f(X_1,\ldots,X_n)=(\delta_{1j_0},\ldots,\delta_{j_0j_0},\ldots,\delta_{nj_0}).
$$
$$
{\rm and}\,\,\,\,\,(\alpha_1,\ldots,\alpha_n)=(1,\ldots,1).
$$
Then conditions (i) and (ii) of Corollary 1.3, with the choice $j(i)=j_0$ are satisfied. Since $f(\mathbb{R}^n)\ne\mathbb{R}^n$ it 
must be that the system (\ref{l5}) has non-zero solutions over $\mathbb{R}$. But in this case the system (\ref{l5}) coincides with the system
above, (\ref{l6}). $\qed $

\begin{theorem}
Let the polynomial mapping $f:\,\mathbb{R}^n\rightarrow\mathbb{R}^n$, be given by $f(X_1,\ldots,X_n)=(p_1(X_1,\ldots,X_n),\ldots,
p_n(X_1,\ldots,X_n))$. Suppose that the determinant $\det J(f)(X_1,\ldots,X_n)$ never vanishes in $\mathbb{R}^n$. 
Then $\forall\,j$, $1\le j\le n$ the system:
\begin{equation}
\label{l7}
\overline{p}_j\overline{\frac{\partial p_j}{\partial X_i}}=0,\,\,\,\,\,i=1,\ldots,n,
\end{equation}
has non-trivial solutions over $\mathbb{R}$.
\end{theorem}
\noindent
{\bf Proof.} \\
Let $j=j_0$ be such that the system (\ref{l7}) has only the zero solution over $\mathbb{R}$. We will arrive at a contradiction by
showing that this assumption implies on the one hand, $f(\mathbb{R}^n)=\mathbb{R}^n$, and it also implies, on the other hand, 
$f(\mathbb{R}^n)\ne\mathbb{R}^n$. \\
1) We first prove that $f(\mathbb{R}^n)=\mathbb{R}^n$. To see that, we use Corollary 1.3 with:
$$
g_{ij}(X_1,\ldots,X_n)=\frac{\partial p_j}{\partial X_i}\,\,\,\,\,{\rm for}\,\,i,j=1,\ldots,n.
$$
We can assume without losing the generality that:
\begin{equation}
\label{l8}
\deg g_{ij_0}=\deg p_{j_0}-1,\,\,\,\,\,i=1,\ldots,n.
\end{equation}
For the assumptions of our theorem as well as the conclusion $f(\mathbb{R}^n)=\mathbb{R}^n$, are invariant with respect to a
real, non-singular change of the variables. More precisely, instead of working with the original mapping, $f(X_1,\ldots,X_n)$, we 
could have, first performed a change of the variables, as follows:
\begin{equation}
\label{l9}
X_j=\sum_{i=1}^n a_{ij}U_i,\,\,\,\,\,j=1,\ldots,n,
\end{equation}
where $(a_{ij})_{i,j=1,\ldots,n}$ is a real non-singular matrix. Then we could have proved that the mapping given by 
$F(U_1,\ldots,U_n)=f(X_1,\ldots,X_n)$ is epimorphic and that would have implied that the original mapping $f(X_1,\ldots,X_n)$ is epimorphic.
The linear transformation we choose in equation (\ref{l9}) is such that $a_{ij}\ne 0$ for all $i,j=1,\ldots,n$. With this choice of the
linear transformation it is clear that generically (in the $a_{ij}\ne 0$), each of the components $\tilde{p}_j(U_1,\ldots,U_n)=p_j(X_1,\ldots,X_n)$,
$j=1,\ldots,n$, of the mapping $F(U_1,\ldots,U_n)$ has the property that for each $i=1,\ldots,n$ it contains all the monomials of
the form $aU_1^{m_1}\ldots U_n^{m_n}$ where $a\ne 0$, and where $\sum_{k=1}^n m_k=\deg p_j$, and $m_i\ne 0$. This justifies
equation (\ref{l8}). Next we choose in Corollary 1.3 the following: For $j\ne j_0$ we take $\alpha_j=1$. We choose the positive integer $\alpha_{j_0}$
so large that $\alpha_{j_0}\deg p_{j_0}+(\deg p_{j_0}-1)$ is strictly larger than $\deg p_j+\deg g_{ij}$ for $i=1,\ldots,n$ and $j\ne j_0$.
Also $\alpha_{j_0}$ is such that $\alpha_{j_0}\deg p_{j_0}+(\deg p_{j_0}-1)$ is an odd integer. That is always possible to do: If $\deg p_{j_0}$
is an even integer, then there is no other restriction on $\alpha_{j_0}$ (except for being large enough). If $\deg p_{j_0}$ is an odd integer,
then $\alpha_{j_0}$ must also be an odd integer. Now conditions (i) and (ii) of Corollary 1.3 are satisfied with $j(i)=j_0$. The system 
(\ref{l5}) Corollary 1.3 reduces to the system (\ref{l7}) with $j=j_0$ and so by part (a) of Corollary 1.3 it follows that $f(\mathbb{R}^n)=\mathbb{R}^n$. \\
2) In order to conclude the proof of Theorem 1.8, we now prove that the existence of such a $j_0$ implies that $f(\mathbb{R}^n)\ne\mathbb{R}^n$.
We may assume that $\overline{p}_{j_0}(X_1,\ldots,X_n)\ge 0$ $\forall\,(X_1,\ldots,X_n)\in\mathbb{R}^n$, and there is an equality 
$\overline{p}_{j_0}(X_1^0,\ldots,X_n^0)=0$ if and only if $(X_1^0,\ldots,X_n^0)=(0,\ldots,0)$. Let us denote $d=\deg\overline{p}_{j_0}$.
We claim that $\forall\,i$, $1\le i\le n$ we have $\deg_{X_i} \overline{p}_{j_0}=d$: \\
For let $\overline{p}_{j_0}(X_1,\ldots,X_n)=\sum_{k=0}^N h_k(X_1,\ldots,\hat{X}_i,\ldots,X_n)X_i^k$ where $h_k$ is an homogeneous polynomial
in $(X_1,\ldots,\hat{X_i},\ldots,X_n)$ of degree $d-k$. Then $\overline{p}_{j_0}(0,\ldots,0,X_i,0,\ldots,0)\equiv 0$ for any choice of $X_i$ 
which is impossible. Hence we obtain:
\begin{equation}
\label{l10}
\overline{p}_{j_0}(X_1,\ldots,X_n)=\sum_{i=1}^n \lambda_iX_i^d+h(X_1,\ldots,X_n),
\end{equation} 
where $\lambda_i>0$, $\forall\,i, 1\le i\le n$ and where $h$ is homogeneous of degree $d$ such that $\deg_{X_{i}}h<d$, $\forall\,i$,
$1\le i\le n$. Since $\overline{p}_{j_0}\ge 0$ it follows that $d$ is an even integer and now equation (\ref{l10}) implies the
existence of an $M>0$ such that $\forall\,(X_1,\ldots,X_n)\in\mathbb{R}^n$ we have $p_{j_0}(X_1,\ldots,X_n)\ge -M$. Hence we
conclude that $f(\mathbb{R}^n)\ne\mathbb{R}^n$. Now the proof of the theorem is completed. $\qed $ \\

\begin{theorem}
Let the polynomial mapping $f:\,\mathbb{R}^n\rightarrow\mathbb{R}^n$, be given by $f(X_1,\ldots,X_n)=(p_1(X_1,\ldots,X_n),\ldots,
p_n(X_1,\ldots,X_n))$. Suppose that the determinant $\det J(f)(X_1,\ldots,X_n)$ never vanishes in $\mathbb{R}^n$. If there is an 
even integral vector $(\alpha_1,\ldots,\alpha_n)\in (2\mathbb{Z}^+)^n$ such that the induced homogeneous system of:
\begin{equation}
\label{l11}
\sum_{j=1}^n\alpha_j\cdot\left(p_j(X_1,\ldots,X_n)\right)^{\alpha_j-1}\frac{\partial p_j}{\partial X_i}=0,\,\,\,\,\,i=1,\ldots,n,
\end{equation} 
has only the zero solution over $\mathbb{R}$, then $f(\mathbb{R}^n)=\mathbb{R}^n$.
\end{theorem}
\noindent
{\bf Proof.} \\
Let us consider the following polynomial: $F(X_1,\ldots,X_n)=\sum_{j=1}^n\left(p_j(X1,\ldots,X_n)\right)^{\alpha_j}$.
Since $p_j(X_1,\ldots,X_n)\in\mathbb{R}[X_1,\ldots,X_n]$, $\forall\,j=1,\ldots,n$ and since the vector $(\alpha_1,\ldots,\alpha_n)$
is an even integral vector, it follows that $\deg F$ is an even integer. Say $\deg F=2N$ for some $N\in\mathbb{Z}^+$. Clearly,
the assumptions as well as the conclusion  of Theorem 1.9 are invariant with respect to a real non-singular linear change of
the variables. Thus, as we explained in the proof of Theorem 1.8 we can assume that:
$$
\deg\left(\frac{\partial F}{\partial X_i}\right)=\deg F -1=2N-1,\,\,\,\,\,i=1,\ldots,n.
$$
Let us take in Theorem 1.2:
$$
g_{ij}(X_1,\ldots,X_n)=\frac{\partial p_j}{\partial X_i},\,\,\,{\rm for}\,\,\,i,j=1,\ldots,n.
$$
The vector of integers in Theorem 1.2 will be $(\alpha_1-1,\ldots,\alpha_n-1)$, and the $j$'th component of the mapping in
Theorem 1.2 will be:
$$
\left(\alpha_j\right)^{1-\alpha_j^{-1}}p_j(X_1,\ldots,X_n).
$$
These satisfy the conditions (i) and (ii) of Theorem 1.2 and now part (a) of Theorem 1.2 implies that $f(\mathbb{R}^n)=\mathbb{R}^n$. $\qed $ \\
\\
{\bf Pinchuk's example.} \\
Pinchuk defined the following: \\
$$
t=xy-1,\,\,\,s=1+xt,\,\,\,h=ts,\,\,\,f=s^2(t^2+y),
$$
and then set,
$$
p=h+f,\,\,\,q=-t^2-6th(h+1)-u(f,h),
$$
where
$$
u=A(h)f+B(h),
$$
$$
A=h+\frac{1}{45}(13+15h)^3,
$$
$$
B=4h^3+6h^2+\frac{1}{2}h^2+\frac{1}{2700}(13+15h)^4.
$$
Thus we have:
$$
\deg h=5,\,\,\,\deg f=10,\,\,\,\deg p=10,\,\,\,\deg q=25.
$$
Pinchuk's example is the following mapping:
$$
(p,q)=(x^6y^4-2x^5y^3+\ldots,\frac{15^3}{45}x^{15}y^{10}+\ldots).
$$
We are interested only in the leading homogeneous components. Thus:
$$
p=x^6y^4+\ldots,\,\,\,\frac{\partial p}{\partial x}=6x^5y^4+\ldots,\,\,\,\frac{\partial p}{\partial y}=4x^6y^3+\ldots.
$$
$$
q=\frac{15^3}{45}x^{15}y^{10}+\ldots,\,\,\,\frac{\partial q}{\partial x}=\frac{15^4}{45}x^{14}y^{10}+\ldots,\,\,\,
\frac{15^3\cdot 10}{45}x^{15}y^9+\ldots.
$$
There are, in this case, two homogeneous systems of equations in (\ref{l7}) of Theorem 1.8:
$$
\overline{p}\frac{\partial\overline{p}}{\partial x}=\overline{p}\frac{\partial\overline{p}}{\partial y}=0,
$$
and
$$
\overline{q}\frac{\partial\overline{q}}{\partial x}=\overline{q}\frac{\partial\overline{q}}{\partial y}=0.
$$
These reduce to:
$$
\begin{array}{ccccc} x^{11}y^8 &=&x^{12}y^{7} &=& 0 \\ x^{29}y^{20} &=& x^{30}y^{19} &=& 0 \end{array}.
$$
Thus both systems have non-zero solutions:
$$
\{(0,y)\,|\,y\in\mathbb{R}\}=\{(x,0)\,|\,x\in\mathbb{R}\},
$$
as should be the case according to Theorem 1.8.

\begin{remark}
The Pinchuk construction gives coordinates with a single element as their highest homogeneous component. This element
has the form $\alpha x^my^k$ where $\alpha\in\mathbb{R}^{\times}$, $m,k\ge 1$. Thus the equations in (\ref{l7}) of 
Theorem 1.8 are of the form:
$$
x^my^k\cdot x^{m-1}y^k=x^my^k\cdot x^my^{k-1}=0,
$$
i.e.
$$
x^{2m-1}y^{2k}=x^{2m}y^{2k-1}=0,
$$
and so the solution set is the union of both axis:
$$
\{(0,y)\,|\,y\in\mathbb{R}\}=\{(x,0)\,|\,x\in\mathbb{R}\},
$$
which, of course, is non-trivial in agreement with Theorem 1.8.
\end{remark}

\noindent
{\it Ronen Peretz \\
Department of Mathematics \\ Ben Gurion University of the Negev \\
Beer-Sheva , 84105 \\ Israel \\ E-mail: ronenp@math.bgu.ac.il} \\ 
 
\end{document}